\documentclass[11pt]{article}
\usepackage{amsfonts}
\usepackage{latexsym}
\usepackage{amsmath}
\newtheorem{theorem}{Theorem}[section]
\newtheorem{lemma}[theorem]{Lemma}

\newtheorem{corollary}[theorem]{Corollary}
\newtheorem{definition}[theorem]{Definition}

\newcommand{\pa}{\hskip 0.55truecm}
\makeatletter
\newcommand{\rank}{\mathop{\operator@font rank}}
\makeatother

\begin{document}
\small

\title{\bf Lower bound for the maximal number of facets of a 0/1 polytope}

\author{D.\ Gatzouras, A.\ Giannopoulos and N.\ Markoulakis}

\date{}

\maketitle

\begin{abstract}
\footnotesize  Let $f_{n-1}(P)$ denote the number of facets of a
polytope $P$ in ${\mathbb R}^n$. We show that there exist 0/1
polytopes $P$ with $$f_{n-1}(P)\geq\left (\frac{cn}{\log^2
n}\right )^{n/2}$$ where $c>0$ is an absolute constant. This
improves earlier work of B\'{a}r\'{a}ny and P\'{o}r on a question
of Fukuda and Ziegler.
\end{abstract}

\section{Introduction}

\pa The aim of this article is to give a lower bound for the
maximal possible number of facets of a 0/1 polytope in ${\mathbb
R}^n$. By definition, a 0/1 polytope is the convex hull of a
subset of the vertices of $[0,1]^n$.

In general, if $P$ is a polytope in ${\mathbb R}^n$, we write
$f_{n-1}(P)$ for the number of its facets. Let $$g(n):=\max\big\{
f_{n-1}(P_n): P_n\;\hbox{a 0/1 polytope in}\;{\mathbb R}^n\big\}
.\leqno (1.1)$$ Fukuda and Ziegler (see \cite{F}, \cite{KRSZ},
\cite{Z}) asked what the behaviour of $g(n)$ is as
$n\rightarrow\infty $. The best known upper bound to date is
$$g(n)\leq 30(n-2)! \leqno (1.2)$$ (for $n$ large enough), which
is established by Fleiner, Kaibel and Rote in \cite{FKR}. We will
study lower bounds.  A major breakthrough in this direction was
made by B\'{a}r\'{a}ny and P\'{o}r in \cite{BaPor}; they proved
that $$g(n)\geq\left (\frac{cn}{\log n}\right )^{n/4},\leqno
(1.3)$$ where $c>0$ is an absolute constant. We will show that the
exponent $n/4$ can in fact be improved to $n/2$:

\begin{theorem} There exists a constant $c>0$ such that
$$g(n)\geq \left (\frac{cn}{\log^2 n}\right )^{n/2}.\leqno (1.4)$$
\end{theorem}

It is interesting to compare this estimate with the known bounds
for the expected number of facets of the convex hull $P_{N,n}$ of
$N$ independent random points which are uniformly distributed on
the sphere $S^{n-1}$. In \cite{BMT} it is shown that there exist
two constants $c_1,c_2>0$, such that $$\left (
c_1\log\frac{N}{n}\right )^{n/2}\leq {\mathbb
E}[f_{n-1}(P_{N,n})]\leq \left ( c_2\log\frac{N}{n}\right
)^{n/2}\leqno (1.5)$$ for all $n$ and $N$ satisfying $2n\leq N\leq
2^n$. In the case of 0/1 polytopes, $N$ can be as large as $2^n$,
therefore one might conjecture that $g(n)$ is of the order of
$n^{n/2}$. Theorem 1.1 gives a lower bound which is ``practically
of this order": for every $\varepsilon >0$ one has $$g(n)
>n^{( 0.5-\varepsilon )n}\leqno (1.6)$$ if $n$ is
large enough.

\smallskip

The existence of 0/1  polytopes with many facets will be
established by a refinement of the probabilistic method developed
in \cite{BaPor}. It will be more convenient to work with $\pm 1$
polytopes (i.e., polytopes whose vertices are sequences of signs).
Let $X_1,\ldots ,X_n$ be independent and identically distributed
$\pm 1$ random variables, defined on some probability space
$(\varOmega,\mathcal{F},\mathbb{P})$, with distribution $$
\mathbb{P}(X=1)=\mathbb{P}(X=-1)=\tfrac{1}{2}.$$  Set
$\vec{X}=(X_1,\ldots ,X_n)$ and, for a fixed $N$ satisfying $n <
N\leq 2^n$, consider $N$ independent copies $\vec{X}_1,\ldots
,\vec{X}_N$ of $\vec{X}$. This procedure defines the random 0/1
polytope $$K_N={\rm conv}\{\vec{X}_1,\ldots ,\vec{X}_N\}.\leqno
(1.7)$$ Note that $K_N$ has at most $N$ vertices (it may happen
that some repetitions occur).

Under some restrictions on the range of values of $N$, we will
obtain a lower bound for the expected number of facets ${\mathbb
E}[f_{n-1}(K_N)]$, for each fixed $N$. In particular, we have:

\begin{theorem} There exist two positive constants $a$ and $b$ such
that: for all sufficiently large $n$, and all $N$ satisfying
$n^a\leq N\leq \exp ( bn/\log n)$, there exists a $0/1$ polytope
$K_N$ in ${\mathbb R}^n$ with $$f_{n-1}(K_N)\geq \left (\frac{\log
N}{a\log n}\right )^{n/2}.\leqno (1.8)$$
\end{theorem}

\noindent It is clear that Theorem 1.1 follows: one only has to
choose $N=\lfloor \exp (bn/\log n)\rfloor $.

\medskip

\section{Preliminaries}

\pa We first fix some standard notation. We work in ${\mathbb
R}^n$ which is equipped with a Euclidean structure $\langle\cdot
,\cdot\rangle $. We denote by $\|\cdot \|_2$ the corresponding
Euclidean norm, by $\|\cdot \|_\infty$ the $\max$-norm, and write
$B_2^n$ for the Euclidean unit ball and $S^{n-1}$ for the unit
sphere. Volume, surface area and the cardinality of a finite set
are denoted by $|\cdot |$ (this will cause no confusion). All
logarithms are natural. Whenever we write $a\simeq b$, we mean
that there exist absolute constants $c_1,c_2>0$ such that
$c_1a\leq b\leq c_2a$. The letters $c,c^{\prime }, C, c_1, c_2$
etc. denote absolute positive constants which may change from line
to line (however, in most places we will try to specify the
absolute constants involved).

As mentioned in the introduction, the proof of Theorem 1.2 will be
modeled on the approach of B\'{a}r\'{a}ny and P\'{o}r. This in
turn has its origin in the work of Dyer, F\"{u}redi and McDiarmid
\cite{DFM}, who proved the following: Let $\kappa =2/\sqrt{e}$ and
consider the random polytope $K_N$ defined in (1.7). For every
$\varepsilon \in (0,1)$, $$\lim_{n\rightarrow\infty }\sup\left\{
2^{-n} {\mathbb E}|K_N|\colon N\leq (\kappa -\varepsilon
)^n\right\}=0\leqno (2.1)$$ and $$\lim_{n\rightarrow\infty
}\inf\left \{ 2^{-n} {\mathbb E}|K_N|\colon N\geq (\kappa
+\varepsilon )^n\right\}=1.\leqno (2.2)$$

In order to determine the threshold $N(n)= (2/\sqrt{e})^n$, they
introduced two families of convex subsets of the cube
$C=[-1,1]^n$. For every $\vec{x}\in C$, set
$$q(\vec{x}):=\inf\big\{ {\rm Prob}\bigl(\vec{X}\in H\bigr)\colon
\vec{x}\in H,\ H\hbox{ a closed halfspace}\bigr\}.\leqno (2.3)$$
If $\beta
>0$ then the $\beta $-center of $C$ is defined by
$$ Q^{\beta }=\{\vec{x}\in C\colon q(\vec{x})\geq\exp (-\beta n)\}
; \leqno (2.4)$$ it is easily checked that $Q^{\beta }$ is a
convex polytope.

Next, consider the function $f:(-1,1)\rightarrow {\mathbb R}$ with
$$f(x)=\tfrac{1}{2}(1+x)\log (1+x)+\tfrac{1}{2}(1-x)\log (1-x)
,\leqno (2.5)$$ extend it to a continuous function on $[-1,1]$ by
setting $f(\pm 1)=\log 2$, and for every $\vec{x}=(x_1,\ldots
,x_n)\in C$ set $$F(\vec{x})=\frac{1}{n}\sum_{i=1}^nf(x_i).\leqno
(2.6)$$ The next lemma was proved in \cite[Section 3]{DFM}.

\begin{lemma}
\label{qleqe-nF}For every $\vec{x}\in (-1,1)^n$ we have
$q(\vec{x})\leq \exp (-nF(\vec{x}))$. $\hfill\Box $\end{lemma}

The second family of subsets of $C$ introduced in \cite{DFM} is as
follows: for every $\beta >0$, set $$F^{\beta }=\{ \vec{x}\in
C:F(\vec{x})\leq \beta\}.\leqno (2.7)$$ Since $f$ is a strictly
convex function on $(-1,1)$, it is clear that $F^{\beta }$ is
convex. Lemma 2.1 and the definition of $Q^{\beta }$ show that if
$\vec{x}\in Q^{\beta }\cap (-1,1)^n$ then $F(\vec{x})\leq\beta $.
In other words, we have the following.

\begin{lemma}
\label{qsubsetf} $Q^{\beta }\cap (-1,1)^n\subseteq F^{\beta }$ for
every $\beta
>0$. $\hfill\Box $
\end{lemma}

Observe that as $\beta\rightarrow\log 2$, both $Q^{\beta }$ and
$F^{\beta }$ approach $C$. The main technical step for the proof
of Theorem 1.2 will be to show that the two families are very
close, in the following sense: for a wide range of $\beta $'s, one
has that $$F^{\beta -\varepsilon }\cap\gamma C\subseteq Q^{\beta
}\cap (-1,1)^n\subseteq F^{\beta } ,\leqno (2.8)$$ where $\gamma
>0$ is a (small) absolute constant and $\varepsilon \leq 3\log
n/n$. The estimate on $\varepsilon $ substantially improves (and
at the same time clarifies) \cite[Lemma 4.3]{BaPor} and should be
viewed as the main technical step in our work. The proof is
presented in the next Section. It is based on two ingredients: (1)
Theorem 3.3, which via large deviations exhibits, for each
$\vec{x}\in\partial(F^\beta)$, the precise rate of decay (with
$n$) of the probability that a  randomly chosen vertex of the unit
cube $C=[-1,1]^n$ lies in the halfspace determined by the
hyperplane tangent to $F^\beta$ at $\vec{x}$ which does not
contain $F^\beta$; and (2) on a result of Montgomery-Smith
\cite{SMS}, giving lower bounds for the tails of the distribution
of Rademacher sums.

\medskip

\section{Comparing $F^{\alpha }$ to $Q^{\alpha }$: the main lemma}

\pa Let $n\in {\mathbb N}$ and let $X_1,X_2,\ldots,X_n$ be
independent and identically distributed $\pm 1$ random variables,
defined on some probability space
$(\varOmega,\mathcal{F},\mathbb{P})$, with distribution
$\mathbb{P}(X=1)=\mathbb{P}(X=-1)=\tfrac{1}{2}$.  For
$t\in\mathbb{R}$, let $$ \varphi(t) := \mathbb{E}
\left[e^{tX}\right] = \cosh(t) \leqno (3.1)$$ be the common moment
generating function of the $X_i$'s, and set $$ \psi(t) := \log
\varphi(t) = \log \cosh(t) .\leqno (3.2)$$  Finally, define
$h:(-1,1)\rightarrow {\mathbb R}$ by
$$h(x):=\tfrac{1}{2}\log\left(\frac{1+x}{1-x}\right) .\leqno
(3.3)$$ It is easily checked that $h$ is strictly convex and
strictly increasing on $[0,1)$. It is also easily checked that
$$f(x) = -\psi(h(x)) + x\, h(x)\quad\mbox{and}\quad f'(x)=h(x)
.\leqno (3.4)$$

Given $x_1,\ldots,x_n$ in $(-1,1)$, set $$t_i := h(x_i)\qquad
(i\leq n) ;\leqno(3.5)$$ in the sequel, $t_i$ and $x_i$ will
always be in this relationship. Observe that
$$x_i=\psi'(t_i)=\tanh(t_i)\qquad (i\leq n).\leqno(3.6)$$ Define a
new probability measure $\mathbb{P}_{x_1,\ldots,x_n}$ on
$(\varOmega,\mathcal{F})$, by $$ \mathbb{P}_{x_1,\ldots,x_n}(A) :=
{\mathbb E} \left[ \mathbf{1}_A \exp\left(\sum\limits_{i=1}^n
t_iX_i\right)\right]\,\prod\limits_{i=1}^n
[\varphi(t_i)]^{-1}\leqno (3.7)$$ for $A\in\mathcal{F}$. The next
lemma is verified by direct computation.

\begin{lemma}
\label{moments} Under $\mathbb{P}_{x_1,\ldots,x_n}$, the random
variables $t_1X_1,\ldots,t_nX_n$ are independent, with mean,
variance and absolute central third moment given by
\begin{eqnarray*} {\mathbb E}_{x_1,\ldots,x_n}[t_iX_i] &=& t_ix_i,\\
{\mathbb E}_{x_1,\ldots,x_n}\left[ t_i^2(X_i-x_i)^2\right] &=&
\frac{t_i^2}{\cosh^2(t_i)} ,\\ {\mathbb E}_{x_1,\ldots,x_n} \left[
|t_i(X_i-x_i)|^3\right] &=&
|t_i|^3\frac{\cosh(2t_i)}{\cosh^4(t_i)} ,\end{eqnarray*}
respectively. $\hfill\Box $
\end{lemma}

Set $$\sigma_n^2 := \sum\limits_{i=1}^n {\mathbb
E}_{x_1,\ldots,x_n}\left[ t_i^2(X_i-x_i)^2\right] =
\sum\limits_{i=1}^n \frac{t_i^2}{\cosh^2(t_i)} \leqno (3.8)$$ and
$$ S_n := \frac{1}{\sigma_n} \sum\limits_{i=1}^n t_i(X_i-x_i)
,\leqno (3.9)$$ and let $F_n\colon\mathbb{R}\rightarrow\mathbb{R}$
be the cumulative distribution function of the random variable
$S_n$ under the probability law $\mathbb{P}_{x_1,\ldots,x_n}$:
$$F_n(x) := \mathbb{P}_{x_1,\ldots,x_n} (S_n\leq x)\qquad
(x\in\mathbb{R}) .\leqno (3.10)$$ Write also $\mu_n$ for the
probability measure on $\mathbb{R}$ defined by $$\mu_n(-\infty,x]
:= F_n(x) \qquad(x\in\mathbb{R}) .\leqno(3.11)$$ Finally, set $$
\rho^{(3)}_n := \sum\limits_{i=1}^n {\mathbb E}_{x_1,\ldots,x_n}
\left[ |t_i(X_i-x_i)|^3\right] = \sum\limits_{i=1}^n
|t_i|^3\frac{\cosh(2t_i)}{\cosh^4(t_i)} .\leqno (3.12)$$ Since
$\cosh (2y)\leq 2\cosh^2(y)$, we have that
$$\frac{\rho_n^{(3)}}{\sigma_n^2}\leq 2\max\limits_{1\leq i\leq n}
|t_i|.\leqno (3.13)$$  Notice also that $${\mathbb
E}_{x_1,\ldots,x_n}[S_n]=0\quad\mbox{and}\quad
\mbox{Var}_{x_1,\ldots,x_n} [S_n]=1 .$$

By (3.7), we have that
\begin{multline*}\mathbb{P}\left(\sum\limits_{i=1}^n t_i(X_i-x_i)\geq 0\right)
\\ = {\mathbb E}_{x_1,\ldots,x_n} \left[
\mathbf{1}_{[0,\infty)}\left(\sum_{i=1}^n t_i(X_i-x_i)\right)\,
\exp\left(-\sum_{i=1}^n t_i X_i\right)\right] \prod\limits_{i=1}^n
\varphi(t_i);\end{multline*} hence, using (3.4), (3.8), (3.9) and
(3.11), we see that $$ \mathbb{P}\left(\sum\limits_{i=1}^n
t_i(X_i-x_i)\geq 0\right) = \int_{[0,\infty)} e^{-\sigma_n u}\,
d\mu_n(u)\, \exp\left(\sum_{i=1}^n [\psi(t_i) - t_ix_i)]\right) ,
$$ which then, upon using (3.4) and (3.5), yields that,
$$\mathbb{P}\left(\sum\limits_{i=1}^n t_i(X_i-x_i)\geq 0\right) =
\int_{[0,\infty)} e^{-\sigma_n u}\, d\mu_n(u)\, \exp\left(-
\sum_{i=1}^n f(x_i)\right) .\leqno(3.14) $$

Next, write $$\phi(x):=\dfrac{1}{\sqrt{2\pi}} e^{-x^2/2}\qquad
(x\in\mathbb{R})\leqno(3.15)$$ for the standard gaussian density,
$$\varPhi(x) := \int_{-\infty}^x \phi(y)\, d y \qquad
(x\in\mathbb{R})\leqno (3.16)$$ for the standard gaussian c.d.f.,
and $\mu$ for the standard gaussian probability measure:
$\mu(-\infty,x] := \varPhi(x)$, $x\in\mathbb{R}$.  By the
Berry--Esseen theorem \cite[Theorem XVI.5.2]{Fe} and Lemma
\ref{moments}, one has $$|F_n(x) - \varPhi(x)|\leq 6\,
\frac{\rho_n^{(3)}}{\sigma_n^3} \leqno (3.17)$$ for all $x$,
whence by (3.13), $$|F_n(x) - \varPhi(x)|\leq
\frac{12}{\sigma_n}\, \max\limits_{1\leq i \leq n} |t_i| \leqno
(3.18)$$ for all $x$.

\begin{lemma}
\label{inequality} The following inequality holds: $$ \left
|\int_{(0,\infty )} e^{-\sigma_nu}\, d\mu_n(u) - \int_{(0,\infty
)} e^{-\sigma_nu}\, d\mu(u)\right | \leq \frac{24}{\sigma_n}\,
\max\limits_{1\leq i \leq n} |t_i| .$$
\end{lemma}

\noindent {\it Proof}. Since \begin{eqnarray*}\int_{(0,\infty )}
e^{-\sigma_nu}\, d\mu_n(u) &=& \int_0^\infty \mu_n\left(\{u\colon
e^{-\sigma_nu}\mathbf{1}_{(0,\infty)}(u)\geq r\}\right)\, dr \\
&=& \int_0^1 \mu_n\left(0,-\sigma_n^{-1}\log r\right]\, dr\\ &=&
\int_0^1 \left[ F_n\left(-\sigma_n^{-1}\log r\right) -
F_n(0)\right]\, dr ,\end{eqnarray*} and similarly
\[\int_0^\infty e^{-\sigma_n u}\, d\mu(u) = \int_0^1 \left[
\varPhi\left(-\sigma_n^{-1}\log r\right) - \varPhi(0)\right]\, dr
,\] the result follows from (3.18). $\hfill\Box $

\medskip

We will use the estimates $$ \frac{1}{x}\, m_1(x)\, e^{-x^2/2}
\leq \int_x^\infty \phi(u)\, du \leq \frac{1}{x}\, m_2(x)\,
e^{-x^2/2}\qquad (x>0) ,$$ where $$m_1(x) =
\dfrac{1}{\sqrt{2\pi}}\,
\frac{2x}{x+\sqrt{x^2+4}}\qquad\mbox{and}\qquad m_2(x) =
\dfrac{1}{\sqrt{2\pi}}\, \frac{4x}{3x+\sqrt{x^2+8}}\leqno (3.19)$$
(see \cite[p.\ 17]{IMK} and, for the upper estimate, \cite{SzW}).
Since $$\int_0^\infty e^{-\sigma_n u}\, d\mu(u) =
\frac{e^{\sigma_n^2/2}}{\sqrt{2\pi}}\int_0^\infty
e^{-(\sigma_n+u)^2/2}\, du = e^{\sigma_n^2/2}
\int_{\sigma_n}^\infty \phi(u)\, du ,$$ it follows that
$$\frac{m_1(\sigma_n)}{\sigma_n} \leq \int_0^\infty e^{-\sigma_n
u}\, d\mu(u) \leq \frac{m_2(\sigma_n)}{\sigma_n} . \leqno (3.20)$$
Combining (3.20) with Lemma \ref{inequality}, we obtain the
estimates $$\dfrac{m_1(\sigma_n)}{\sigma_n} -
\dfrac{24}{\sigma_n}\, \max\limits_{1\leq i\leq n}|t_i| \leq
\int_{(0,\infty )}e^{-\sigma_n u}\, d\mu_n(u) \leq
\dfrac{m_2(\sigma_n)}{\sigma_n} + \dfrac{24}{\sigma_n}\,
\max\limits_{1\leq i\leq n}|t_i| . \leqno (3.21)$$  Equation
(3.14) then yields the following estimates:

\begin{theorem}
\label{result} Let $x_1,\ldots,x_n\in (-1,1)$ and $t_i=h(x_i)$,
$i=1,\ldots,n$. Then, $$ \mathbb{P}\left(\sum\limits_{i=1}^n
t_i(X_i-x_i)\geq 0\right) \leq \dfrac{1}{\sigma_n}\, e^{-n
F(\vec{x})} \left(m_2(\sigma_n) + 48 \max\limits_{1\leq i\leq n}
|t_i|\right) \leqno(3.22) $$ and
$$\mathbb{P}\left(\sum\limits_{i=1}^n t_i(X_i-x_i)\geq 0\right)
\geq \dfrac{1}{\sigma_n}\, e^{-n F(\vec{x})} \left(m_1(\sigma_n) -
24 \max\limits_{1\leq i\leq n} |t_i|\right) . \leqno(3.23) $$
\end{theorem}

\noindent {\it Proof}. Observe that the first factor on the right
hand-side of (3.14) is $$\int_{[0,\infty)} e^{-\sigma_n u}\,
d\mu_n(u) \geq \int_{(0,\infty)} e^{-\sigma_n u}\, d\mu_n(u). $$
Hence, we obtain the second inequality (3.23) by combining (3.14)
with (3.21). For the first inequality, first observe that $$
\int_{[0,\infty)} e^{-\sigma_n u}\, d\mu_n(u) = \int_{(0,\infty)}
e^{-\sigma_n u}\, d\mu_n(u) + \mathbb{P}_{x_1,\ldots,x_n} (S_n =
0) ,$$ and then use (3.18) to obtain, $$
\mathbb{P}_{x_1,\ldots,x_n} (S_n=0) \leq F_n(\epsilon) -
F_n(-\epsilon)\leq  \varPhi(\epsilon) - \varPhi(-\epsilon) +
\frac{24}{\sigma_n}\, \max\limits_{1\leq i \leq n} |t_i| $$ for
all $\epsilon >0$. Thus $$ \mathbb{P}_{x_1,\ldots,x_n} (S_n=0)
\leq \frac{24}{\sigma_n}\, \max\limits_{1\leq i \leq n} |t_i|
,\leqno(3.24)$$ and the first inequality (3.23) follows now as
well, using (3.14), (3.21) and (3.24) this time. $\hfill\Box$

\begin{corollary}
\label{corollary} Let $\delta\in (0,1)$. If $x_1,\ldots ,x_n\in
(-\delta ,\delta )$ and $t_i=h(x_i)$, $i=1,\ldots ,n$, then
\begin{multline*} \mathbb{P}\left(\sum\limits_{i=1}^n t_i(X_i-x_i)\geq
0\right) \leq \frac{1}{\sqrt{n F(\vec{x})}}\, e^{-nF(\vec{x}) }\\
\times \frac{\cosh(h(\delta))}{\sqrt{2}} \left( m_2\! \left
(h(\delta) (f(\delta))^{-1/2} \sqrt{n F(\vec{x})}\right ) + 48\,
h(\delta )\right ) .\end{multline*} and
\begin{multline*} \mathbb{P}\left(\sum\limits_{i=1}^n t_i(X_i-x_i)\geq
0\right) \geq \frac{1}{\sqrt{n F(\vec{x})}}\, e^{-nF(\vec{x}) }\\
\times \frac{\sqrt{f(\delta )}}{h(\delta)} \left( m_1\! \left
((\cosh(h(\delta)))^{-1} \sqrt{2 n F(\vec{x})}\right ) - 24\,
h(\delta)\right ) .\end{multline*}
\end{corollary}

\noindent {\it Proof}. Recall (3.8). The function
$$g(t)=\frac{f(\tanh (t))}{t^2}=-\frac{1}{t^2}\log\cosh
(t)+\frac{\tanh (t)}{t}\leqno (3.25)$$ is strictly decreasing on
$[0,\infty )$ and $\lim_{t\rightarrow 0}g(t)=\tfrac{1}{2}$ (see
\cite[Lemma 6.1]{BaPor}). It follows that $$\frac{f(\delta
)}{h^2(\delta )}t_i^2\leq f(x_i)\leq\frac{1}{2}t_i^2\leqno
(3.26)$$ for all $i\leq n$, since also $h$ is increasing on
$[0,1)$. Since $1\leq \cosh^2(t_i)\leq\cosh^2(h(\delta ))$, we get
that $$\frac{f(\delta )}{h^2(\delta )}\, \sigma_n^2\leq n
F(\vec{x})\leq \frac{\cosh^2( h(\delta))}{2}\, \sigma_n^2.\leqno
(3.27)$$ Finally, we also have that $$\max_{1\leq i\leq
n}|t_i|\leq h(\delta ) .\leqno (3.28)$$  Inserting these estimates
into the estimates of Theorem 3.3, and using the fact that $m_1$
and $m_2$ are increasing on $[0,\infty )$, concludes the proof.
$\hfill\Box
$

\medskip

The upper bound (3.22), and its equivalent in Corollary
\ref{corollary}, will not be used in the sequel and are only given
for completeness. Notice, however, that these bounds subsume Lemma
2.1.

\begin{corollary}
There exist $\gamma\in (0,1)$ and $k=k(\gamma )\in {\mathbb N}$
with the following property: For every $n\in {\mathbb N}$, if
$x_1,\ldots ,x_n\in (-\gamma ,\gamma )$ are such that
$\sum_{i=1}^nf(x_i)\geq k(\gamma )$, and if $t_i=h(x_i)$,
$i=1,\ldots ,n$, then $$ \mathbb{P}\left(\sum\limits_{i=1}^n
t_i(X_i-x_i)\geq 0\right) \geq \frac{\sqrt{f(\gamma )}}{10\,
h(\gamma )}\times\frac{1}{\sqrt{nF(\vec{x})}}\,
e^{-nF(\vec{x})}.\leqno(3.29)$$
\end{corollary}

\noindent {\it Proof}. First choose $\gamma \in (0,1)$ so that
$24h(\gamma )\leq (2\sqrt{2\pi })^{-1}$; this is possible because
$\lim_{\delta\rightarrow 0}h(\delta )=0$.

We know that $m_1$ increases to $(2\pi)^{-1/2}$ as
$x\rightarrow\infty $; so, there exists $k=k(\gamma )\in {\mathbb
N}$ such that $$m_1\left(\frac{\sqrt{2k(\gamma )}}{\cosh \big (
h(\gamma )\big )}\right)\geq \frac{5}{6\sqrt{2\pi }}.$$ From the
second assertion of Corollary 3.4 we easily check (3.29) for all
$x_i\in (-\gamma ,\gamma )$ with
$nF(\vec{x})=\sum_{i=1}^nf(x_i)\geq k(\gamma )$. $\hfill\Box $

\medskip

We shall also make essential use of the following result of
Montgomery-Smith from \cite{SMS}:

\begin{lemma}
\label{MontSm} There exists a universal constant $c$, such that,
for all $n\in\mathbb{N}$ and any $s_1,\ldots,s_n\in\mathbb{R}$,
the inequality $$ \mathbb{P}\left(\sum_{i=1}^n s_i X_i \geq c^{-1}
a \|\vec{s}\|_2\right) \geq c^{-1} e^{-ca^2} $$ holds for all
$a>0$ with $a\leq \|\vec{s}\|_2/\|\vec{s}\|_\infty$, where
$\vec{s}=(s_1,\ldots,s_n)$. $\hfill\Box
$
\end{lemma}

\begin{definition} {\rm In the sequel we fix a constant $\gamma\in (0,1)$ which
satisfies Corollary 3.5 and also $\gamma\leq \tanh (c^{-1})$,
where $c$ is the universal constant of Lemma \ref{MontSm}. For
example, since one can have $c=4\log 12$ in Lemma 3.6, we may
take}
$$\gamma = \min\left\{\tanh\left(\frac{1}{48\sqrt{2\pi }}\right) ,\tanh
\left (\frac{1}{4\log 12}\right )\right\}
=\tanh\left(\frac{1}{48\sqrt{2\pi }}\right).\leqno (3.30)$$
\end{definition}

\noindent We shall also use the following lemma (see \cite[Lemma
8.2]{BaPor}).

\begin{lemma}Let $\gamma\in (0,1)$ and $s_i>0$, $i=1,\ldots,m$.
Then, for $m>2/(1-\gamma)$, \begin{eqnarray*} \mathbb{P} \left (
\sum_{i=1}^m s_i(X_i-\gamma )\geq 0\right ) &\geq&
\sqrt{\frac{2}{\pi(1-\gamma^2)}}\,
\frac{1-\gamma-2m^{-1}}{1+\gamma +2m^{-1}}\, \\ & &\times\,
\exp\left(\frac{1}{12m+1}-\frac{1}{12m}\, \frac{4}{1-(\gamma
+2m^{-1})^2}\right)\, \\ & &\times\, \frac{1}{m^{3/2}}\, \exp (-m
f(\gamma )).
\end{eqnarray*}
In particular, for $\gamma\leq\tanh \big ((48\sqrt{2\pi})^{-1}\big
)$ and $m>2$, $$ \mathbb{P} \left ( \sum_{i=1}^m s_i(X_i-\gamma
)\geq 0\right ) \geq c(\gamma)\, m^{-3/2}\, e^{-mf(\gamma)} ,$$
where $c(\gamma)>0$ is given by $$c(\gamma) =
\sqrt{\frac{2}{\pi(1-\gamma^2)}}\, \frac{1-3\gamma}{5+3\gamma}\,
\exp(-[9-(3\gamma+2)^2]^{-1}) .$$
\end{lemma}

\noindent {\it Proof}. This is proved as Lemma 8.2 of
\cite{BaPor}. The argument there gives $$\mathbb{P} \left (
\sum_{i=1}^m s_i(X_i-\gamma )\geq 0\right ) \geq \frac{1}{m}\,
\frac{1}{2^m}\, \sum\limits_{\frac{\gamma+1}{2}\leq \frac{k}{m}
\leq 1} {m\choose k} \geq \frac{1}{m}\, \frac{1}{2^m}\, {m\choose
k_\gamma} , $$ where $k_\gamma := \lceil m (\gamma+1)/2 \rceil$ is
the least integer $\geq m (\gamma+1)/2$; in particular,
$$\frac{k_\gamma}{m}<\frac{1+\gamma+2m^{-1}}{2}.$$ Now using
H.~E.\ Robbins' fine form of the Stirling approximation given in
\cite[II, (9.15)]{Fe1}, and the fact that $$f(\gamma+2m^{-1})\leq
f(\gamma) + 2 m^{-1} h(\gamma+2m^{-1})$$ (by the mean value
theorem and the monotonicity of $f'=h$), yields the result.
$\hfill\Box
$

\begin{theorem} There exists $n_0=n_0(\gamma )\in {\mathbb N}$ with the
following property: If $n\geq n_0$ and $4\log
n/n\leq\alpha\leq\log 2$, then $$ F^{\alpha -\varepsilon_2}\cap
\gamma C\subseteq Q^{\alpha }$$ for some $\varepsilon_2\leq 3\log
n/n$.
\end{theorem}

\noindent {\it Proof}. Fix $\varepsilon_2=3\log n/n$. We need to
check that $q(\vec{x})\geq\exp (-\alpha n)$ for every $\vec{x}$ in
$F^{\alpha -\varepsilon_2}\cap \gamma C$. It suffices to prove
that $$ \mathbb{P}\big (\vec{X}\in H\big )\geq \exp (-\alpha
n)\leqno (3.31)$$ for every halfspace $H$ touching $F^{\alpha
-\varepsilon_2}\cap \gamma C$. In the proof of \cite[Lemma
4.3]{BaPor} it is explained that, for any such halfspace $H$,
there exists $\vec{x}$ on the bounding hyperplane of $H$ such that
$F(\vec{x})=\alpha -\varepsilon_2$. By symmetry we may assume that
$\vec{x}=(x_1,\ldots ,x_n)$ with $0 < x_1\leq x_2\leq \cdots\leq
x_n$.

There exists $n_1\in\{ 1,\ldots,n\}$ such that $x_{n_1}<\gamma $
and $x_{n_1+1}=\gamma $. As in \cite[Lemma 4.3]{BaPor} we set
$\vec{t}=(t_1,\ldots ,t_n)=n\, \nabla F(\vec{x})$ and write
$\vec{t}_{\ast }=(t^{\ast }_1,\ldots ,t^{\ast }_n)$ for the normal
to the bounding hyperplane of $H$. We may assume that
$\vec{t}_{\ast }$ is in the relative interior of the normal cone
to $F^{\alpha -\varepsilon_2}\cap \gamma C$ at $\vec{x}$, whence
$$t_i^{\ast }=t_i=f^{\prime }(x_i) \quad\mbox{if}\quad i\leq
n_1\quad\mbox{and}\quad t_i^{\ast }>t_i=f^{\prime }(\gamma
)\quad\mbox{if}\quad i>n_1.$$

Write
\begin{multline*} \mathbb{P}\left (\sum_{i=1}^nt_i^{\ast }(X_i-x_i)\geq 0\right
) = \mathbb{P}\left (\sum_{i=1}^{n_1}t_i(X_i-x_i) +
\sum_{i=n_1+1}^nt_i^{\ast }(X_i-\gamma)\geq 0\right )\\  \geq
\mathbb{P}\left (\sum_{i=1}^{n_1}t_i(X_i-x_i)\geq 0\right ) \
\mathbb{P}\left (\sum_{i=n_1+1}^nt_i^{\ast }(X_i-\gamma)\geq
0\right ) .\end{multline*}  We estimate the second probability in
the last product using Lemma 3.8: $$\mathbb{P}\left
(\sum_{i=n_1+1}^nt_i^{\ast }(X_i-\gamma)\geq 0\right )\geq
\exp\left( -(n-n_1) f(\gamma) - \tfrac{3}{2} \log(n-n_1) -
c_1(\gamma)\right) .\leqno(3.32)$$ To estimate the first
probability we distinguish two cases:

\medskip

\noindent {\sc Case 1}: $\sum_{i=1}^{n_1} f(x_i)\geq k(\gamma)$.
We may then use Corollary 3.5 to estimate the first probability:
$$ \mathbb{P}\left (\sum_{i=1}^{n_1}t_i(X_i-x_i)\geq 0\right )\geq
\exp\left( - \sum_{i=1}^{n_1} f(x_i) - \tfrac{1}{2} \log
\sum_{i=1}^{n_1} f(x_i) - c_2(\gamma) \right).\leqno(3.33)$$
Combining (3.32) and (3.33) we obtain
\begin{eqnarray*}
\mathbb{P}\left (\sum_{i=1}^nt_i^{\ast }(X_i-x_i)\geq 0\right ) &
\geq &  \exp\left( -\sum_{i=1}^{n_1} f(x_i) - \tfrac{1}{2} \log
\sum_{i=1}^{n_1} f(x_i)\right)\\ & & \times\, \exp\left(-(n-n_1)
f(\gamma) - \tfrac{3}{2} \log(n-n_1) - c(\gamma)\right)\\ &\geq &
\exp\left( -\sum_{i=1}^n f(x_i) - 2 \log n - c(\gamma)\right) \\ &
= & \exp\left( -(\alpha-\varepsilon_2) n - 2 \log n -
c(\gamma)\right) \\ &\geq & \exp(\alpha n ),
\end{eqnarray*}
provided $n$ is large enough to have that $\log n \geq c(\gamma)$.

\medskip

\noindent {\sc Case 2}: $\sum_{i=1}^{n_1} f(x_i) < k(\gamma)$. In
this case we use Lemma \ref{MontSm} to estimate the first
probability.  We have that $$\mathbb{P}\left(
\sum\limits_{i=1}^{n_1} t_i(X_i-x_i)\geq 0\right) =
\mathbb{P}\left( \sum\limits_{i=1}^{n_1} t_iX_i\geq c^{-1} a\,
\sqrt{\sum\limits_{i=1}^{n_1} t_i^2} \right) ,\leqno(3.34)$$ with
$$a=c\, \dfrac{\sum_{i=1}^{n_1} t_ix_i}{\sqrt{\sum_{i=1}^{n_1}
t_i^2}} .$$ To use Lemma \ref{MontSm} we need to check that
$$a\leq\frac{\sqrt{\sum_{i=1}^{n_1} t_i^2}}{\max_{1\leq i\leq n_1} t_i}.$$

The function $h$ is convex on $[0,1)$ and its derivative at $x=0$
is equal to $1$; hence $x\leq h(x)$ for all $x\in [0,1)$.  It
follows that $\sum_{i=1}^{n_1} t_i x_i \leq \sum_{i=1}^{n_1}
t_i^2$, by (3.5).  Therefore, $$ a\leq c\,
\sqrt{\sum\limits_{i=1}^{n_1} t_i^2} \leq c\, h(\gamma)\,
\dfrac{\sqrt{\sum_{i=1}^{n_1} t_i^2}}{\max\limits_{1\leq i\leq
n_1} t_i} ,\leqno(3.35)$$ since also $t_i=h(x_i)\leq h(\gamma)$
for all $i$, by the monotonicity of $h$.  By (3.30) and the
monotonicity of $h$, $c\, h(\gamma)\leq 1$ and $a$ satisfies the
condition of Lemma \ref{MontSm}.  Lemma \ref{MontSm} and (3.34)
then yield the bound $$\mathbb{P}\left( \sum\limits_{i=1}^{n_1}
t_i(X_i-x_i)\geq 0\right)  \geq c^{-1} e^{- c a^2 } \geq
\dfrac{1}{c}\exp\left( - c^3\, \dfrac{h^2(\gamma)}{f(\gamma)}\,
\sum\limits_{i=1}^{n_1} f(x_i) \right), \leqno(3.36)$$ the last
inequality by the first inequality in (3.35) and (3.26). Then,
(3.36) and (3.32) yield the bound
\begin{eqnarray*} & &\mathbb{P}\left( \sum\limits_{i=1}^n
t_i(X_i-x_i)\geq 0\right)\\ & &\quad\geq \exp\left( -(n-n_1)
f(\gamma) - \tfrac{3}{2} \log(n-n_1) - c_1(\gamma)
-c_3(\gamma)\sum\limits_{i=1}^{n_1} f(x_i) - c_4(\gamma)\right)\\
& &\quad \geq \exp\left( - \sum\limits_{i=1}^{n} f(x_i) -
\tfrac{3}{2} \log(n-n_1) - c_1(\gamma) -|1-c_3(\gamma)| k(\gamma)
- c_4(\gamma) \right)\\ & &\quad = \exp\left( -(\alpha -
\varepsilon_2) n - \tfrac{3}{2} \log(n-n_1) - C(\gamma) \right) \\
& &\quad \geq \exp\left( -\alpha n \right) ,
\end{eqnarray*}
provided again that $n$ is large enough to have $\log n\geq
C(\gamma)$.

\medskip

Therefore, in both Cases we have that $$ \mathbb{P}\left
(\sum_{i=1}^nt_i^{\ast }(X_i-x_i)\geq 0\right )\geq\exp \big (
-\alpha n\big ) $$ for $n\geq n_0(\gamma)$, and this completes the
proof. $\hfill\Box
$

\medskip

\section{Weakly sandwiching $K_N$}

\pa The families $\{Q^{\beta }\}$ and $\{F^{\beta }\}$ are related
to the behaviour of the random polytope $K_N$. Fix $N$ with $n <
N\leq 2^n$ and define $\alpha $ by the equation $N=e^{\alpha n}$.
In other words, $$\alpha =\frac{\log N}{n}.\leqno (4.1)$$ In
\cite[Lemma 4.2]{BaPor} (see also \cite[Lemma 2.1(b)]{DFM}), it is
proved that, for some small $\varepsilon_1 (N,n)\in (0,1)$, the
probability ${\rm Prob}\big (Q^{\alpha -\varepsilon_1
}\not\subseteq K_N\big )$ is very small if $n$ is large enough. On
the other hand, in \cite[Lemma 4.4]{BaPor} it is proved that, for
some small $\varepsilon_3 (N,n)\in (0,1)$, at least half of the
surface area of $F^{\alpha +\varepsilon_3 }$ lying in
$\tfrac{1}{10}C$ is missed by the typical $K_N$. In view of
Theorem 3.9 (or \cite[Lemma 4.3]{BaPor}), this means that
$K_N\cap\tfrac{1}{10}C$ is ``weakly sandwiched" between $F^{\alpha
-(\varepsilon_1+\varepsilon_2) }\cap\tfrac{1}{10}C$ and $F^{\alpha
+\varepsilon_3}$. In this Section we provide new estimates for the
parameters $\varepsilon_1$ and $\varepsilon_3$ in the two
statements above.

\begin{lemma} Assume that $e^2n\leq N\leq 2^n$. If $n$ is sufficiently
large, there exists $\varepsilon_1\leq  3\log n/n$ such that $$
{\rm Prob}\big ( K_N\supseteq Q^{\alpha -\varepsilon_1}\big
)>1-2^{-(n-1)}. \leqno (4.2)$$
\end{lemma}

\noindent {\it Proof}. Let $0<\beta <\alpha $. The argument in
\cite[Lemma 2.1(b)]{DFM} gives $$1-{\rm Prob}\big ( K_N\supseteq
Q^{\beta }\big )\leq {N\choose n}\, 2^{-(N-n)}+2\, {N\choose
n}\big (1-e^{-\beta n}\big )^{N-n}.\leqno (4.3)$$ We will bound
both terms on the right handside by $2^{-n}$. Since ${N\choose
n}\leq (eN/n)^n$, we only need to check that $$\left
(\frac{eN}{n}\right )^n2^{-(N-n)}<2^{-n} \leqno (4.4)$$ in order
to bound the first term. This is equivalent to $$ (N-2n)\log
2>n\log \left (\frac{eN}{n}\right ).\leqno (4.5)$$ Since $eN/n\geq
e^3$ and the function $g(x)=\log x/x$ is decreasing for $x\geq e$,
we have that $$ \log \left (\frac{eN}{n}\right )\leq
\frac{3eN}{e^3n}.\leqno (4.6)$$ So, it suffices to check that
$$3e^{-2}N<(\log 2) N-2(\log 2) n.\leqno (4.7)$$ But, $e^2\log
2>2\log 2+3$ and  $N\geq e^2n$; therefore $$ (\log 2
-3e^{-2})N\geq (e^2\log 2-3)n>2(\log 2)n,\leqno (4.8)$$ which
proves (4.4).

\medskip

\noindent Next, we will show that $$ 2\, {N\choose n}\big
(1-e^{-\beta n}\big )^{N-n}<2^{-n} ,\leqno (4.9)$$ provided that
$n$ is large enough and $\alpha -\beta\geq 3\log n/n$. Since
$1-x\leq e^{-x}$, it suffices to check that $$ \left
(\frac{4eN}{n}\right )^n\exp \left ( -e^{-\beta n}(N-n)\right
)<1.\leqno (4.10)$$ Write $\beta =\alpha -\varepsilon $, where
$\varepsilon >0$. Then, $e^{-\beta n}=e^{\varepsilon n}/N$. Since
$n\log (4eN/n)\leq n^2$ (assume that $n\geq 12$) and
$(N-n)/N\geq\tfrac{1}{2}$, we want $$ 2n^2\leq e^{\varepsilon
n}.\leqno (4.11)$$ This is satisfied if $\varepsilon \geq 3\log
n/n$.

From (4.4) and (4.9) we have the lemma. $\hfill\Box $

\begin{lemma} Let $n$ be large enough and assume that $\alpha
<\log 2-12 n^{-1}$. There exists $\varepsilon_3\leq  6/n$ such
that $$ {\rm Prob}\big ( \;|\partial (F^{\alpha
+\varepsilon_3})\cap \gamma C\cap K_N|\geq\tfrac{1}{2}|\partial
(F^{\alpha +\varepsilon_3})\cap \gamma C |\;\big )\leq
\tfrac{1}{100}.\leqno (4.12)$$
\end{lemma}

\noindent {\it Proof}. Let $\beta >\alpha $ and write $\beta
=\alpha +\varepsilon_3$ for some $\varepsilon_3>0$. Let $\vec{x}$
be on the boundary of $F^{\beta }$. If $H$ is a halfspace
containing $\vec{x}$, and if $\vec{x}\in K_N={\rm conv}\{
\vec{X}_1,\ldots ,\vec{X}_N\}$, then there exists $i\leq  N$ such
that $\vec{X}_i\in H$ (otherwise we would have $\vec{x}\in
K_N\subseteq F$, where $F$ is the complementary halfspace). We
write
\begin{eqnarray*}
{\rm Prob} \big (\vec{x}\in K_N\big ) &=& {\rm Prob}\big (\;
\vec{x}\in {\rm conv}\{ \vec{X}_1,\ldots ,\vec{X}_N\}\;\big )\\
&\leq  & {\rm Prob}\big (\; \vec{X}_i\in H \mbox{ for some } 1\leq
i\leq N \;\big )\\ &\leq &\sum_{i=1}^N{\rm Prob}\big (\;
\vec{X}_i\in H\;\big )\\ &=& N\, {\rm Prob}\big ( \vec{X}\in H\big
).
\end{eqnarray*}
Since $H\ni \vec{x}$ was arbitrary, $${\rm Prob}\big (\vec{x}\in
K_N\big )\leq  N \inf\! \big\{ {\rm Prob}\big (\vec{X}\in H\big )
\colon H\ni \vec{x}\big\}=Nq(\vec{x}).\leqno (4.13)$$ From Lemma
2.1 we have $$ Nq(\vec{x})\leq  N\exp (-nF(\vec{x}))=N\exp
(-\alpha n-\varepsilon_3n)=\exp
(-\varepsilon_3n)<\tfrac{1}{200}\leqno (4.14)$$ if
$\varepsilon_3\geq 6/n$. Now $$ {\mathbb E}\; |(\partial (F^{\beta
})\cap \gamma C \cap K_N| \leq \int_{\partial (F^{\beta })\cap
\gamma C }{\rm Prob}\big (\vec{x}\in K_N\big )\, d\vec{x}\leq
\tfrac{1}{200}|\partial (F^{\beta })\cap \gamma C |.\leqno
(4.15)$$ Therefore $$ {\rm Prob}\big (\;|(\partial (F^{\beta
})\cap \gamma C \cap K_N|\geq\tfrac{1}{2}|\partial (F^{\beta
})\cap \gamma C |\;\big )\leq 10^{-2} , \leqno (4.16)$$ by
Markov's inequality. $\hfill\Box $

\medskip

\section{Proof of the Theorem}

We follow the idea of B\'{a}r\'{a}ny and P\'{o}r. We need the
following two facts from \cite{BaPor}:

\begin{lemma} For every $\gamma\in (0,1)$, there exists a constant $c(\gamma
)>0$, such that if $n$ is large enough and $\beta\leq c(\gamma
)/\log n$, then $$ |\partial (F^{\beta })\cap\gamma C|\geq
\tfrac{1}{2} \big(1-\gamma^2\big)^{n-1} (2\beta n)^{(n-1)/2}
|S^{n-1}|.\leqno (5.1)$$
\end{lemma}

\noindent {\it Proof}. We sketch the argument from \cite[Lemma
5.1]{BaPor} in order to make the necessary adjustments. We first
estimate the product curvature $\kappa (\vec{x})$ of the surface
$F(\vec{x})=\beta $ (the formula which appears in \cite{BaPor} is
not exact; we would like to express our gratitude to I.
B\'{a}r\'{a}ny, V. Kaibel and R. Mechtel for kindly pointing out
this point). Let $\nu (\vec{x})=\nabla F(\vec{x})/\|\nabla
F(\vec{x})\|_2$ be the outward unit normal vector of $F^{\beta }$
at $\vec{x}$. Following \cite[Section 2.5]{Schn}, we write
$T_{\vec{x}}F^{\beta}$ for the tangent space of $F^{\beta }$ at
$\vec{x}$, and consider the Weingarten map $W_{\vec{x}} :
T_{\vec{x}} F^{\beta}\rightarrow T_{\vec{x}} F^{\beta}$. This is
the restriction to $T_{\vec{x}} F^{\beta}$ of the differential
$D_{\vec{x}}$ of the map $\vec{x}\mapsto \nu (\vec{x})$. Then
$W_{\vec{x}}$ is symmetric and positive definite, therefore
$$\kappa (\vec{x})=\det W_{\vec{x}}\leq \left (\frac{{\rm
trace}(W_{\vec{x}})}{n-1}\right )^{n-1}\leqno (5.2)$$ by the
arithmetic-geometric means inequality. Let $(a_{ij})_{i,j=1}^n$
denote the matrix of $D_{\vec{x}}$ with respect to the standard
basis of ${\mathbb R}^n$.  As is well known, and readily verified
by direct computation, $\nu(\vec{x})$ is an eigenvector of the
adjoint of $D_{\vec{x}}$, with corresponding eigenvalue $0$; it
follows from this and the fact that the eigenvalues of
$W_{\vec{x}}$ are also eigenvalues of $D_{\vec{x}}$ and none of
them is zero, that ${\rm trace}(W_{\vec{x}})={\rm
trace}(D_{\vec{x}})$. A simple calculation also shows that
$$a_{ii}=\frac{f^{\prime\prime }(x_i)\big (\|n \nabla
F(\vec{x})\|_2^2-(f^{\prime }(x_i))^2\big )}{\|n \nabla
F(\vec{x})\|_2^3}=\frac{h^{\prime }(x_i)\big
(\|\;\vec{t}\;\|_2^2-(h(x_i))^2\big )}{\|\;\vec{t}\;\|_2^3}.\leqno
(5.3)$$ It follows that, if $\vec{x}\in
\partial (F^{\beta })\cap\gamma C$, then
\begin{eqnarray*}
\frac{{\rm trace}(W_{\vec{x}})}{n-1}= \frac{{\rm
trace}(D_{\vec{x}})}{n-1} &=& \sum_{i=1}^n\frac{h^{\prime
}(x_i)\big (\|\;\vec{t}\;\|_2^2-(h(x_i))^2\big
)}{(n-1)\|\;\vec{t}\;\|_2^3}\\ &\leq & h^{\prime }(\gamma
)\frac{n\|\;\vec{t}\;\|_2^2- \sum_{i=1}^n
t_i^2}{(n-1)\|\;\vec{t}\;\|_2^3}\\ &=&\frac{h^{\prime }(\gamma
)}{\|\;\vec{t}\;\|_2},
\end{eqnarray*}
and (5.2) shows that $$\frac{1}{\kappa (\vec{x})}\geq \big\|\;
\vec{t}\; \big\|_2^{n-1}(1-\gamma^2)^{n-1}.\leqno (5.4)$$ Since
also $2f(x_i)\leq t_i^2$, by (3.26), we finally have that
$$\frac{1}{\kappa (\vec{x})}\geq \big(1-\gamma^2\big)^{n-1}
(2n\beta )^{(n-1)/2}.\leqno (5.5)$$

For every $\vec{\theta }\in S^{n-1}$ we write $\vec{x}(\vec{\theta
},\beta )$ for the point on the boundary of $F^{\beta }$ for which
$\vec{t}(\vec{\theta },\beta )=n\, \nabla F(\vec{x}(\vec{\theta
},\beta ))$ is a positive multiple of $\vec{\theta }$. This point
is well-defined and unique if $0<\beta < | {\rm supp} (\vec{\theta
})|\, (\log 2) / n$ (see \cite[Lemma 6.2]{BaPor}). The argument
given in \cite[Lemma 6.3]{BaPor} shows that if $$M =\left\{
\vec{\theta }\in S^{n-1}:\sqrt{\frac{n}{3\log n}}\ \vec{\theta
}\in C\right\}\leqno (5.6)$$ and if $\beta < c\gamma^2/\log n$,
then for every $\vec{\theta }\in M$ we have $\vec{x}(\vec{\theta
},\beta )\in\gamma C$  (observe that if $\vec{\theta }\in M$, we
also have $|{\rm supp}(\vec{\theta })|\geq n/(3\log n)$, and
therefore $\beta < |{\rm supp}(\vec{\theta })|\, (\log 2) / n $ if
$\beta < c\gamma^2/\log n$). A standard computation on the area of
spherical caps shows that $$|M |\geq \tfrac{1}{2}|S^{n-1}|.\leqno
(5.7)$$ Then, we can write $$|\partial (F^{\beta })\cap\gamma C| =
\int_{S^{n-1}}\frac{1}{\kappa (\vec{x})}\, d\vec{\theta } \geq
\int_{M }\frac{1}{\kappa (\vec{x})}\, d\vec{\theta }\geq
\tfrac{1}{2}\big(1-\gamma^2\big)^{n-1}(2n\beta
)^{(n-1)/2}|S^{n-1}|\leqno (5.8)$$ as claimed. $\hfill\Box $

\begin{lemma} Let $\gamma\in (0,0.1)$ and assume that $\beta +\varepsilon <\log 2$.
If $H$ is a halfspace which is disjoint from $F^{\beta }\cap\gamma
C$, then $$ |\partial (F^{\beta +\varepsilon })\cap\gamma C\cap
H|\leq (3\varepsilon n)^{(n-1)/2}|S^{n-1}|.\leqno (5.9)$$
\end{lemma}

\noindent {\it Proof}. Completely similar to the one in
\cite[Lemma 5.2]{BaPor}. $\hfill\Box $

\medskip

\noindent {\bf Proof of Theorem 1.2}: Assume that $n$ is large
enough and set $b=c(\gamma )$, where $c(\gamma )>0$ is the
constant in Lemma 5.1.

Given $N$ with $n^8\leq N\leq\exp \big (bn/\log n\big )$, let
$\alpha =\log N/n$ and let $\gamma \in (0,1)$ be the constant in
(3.30). From Lemma 4.1 there exists $\varepsilon_1\leq 3\log n/n$
such that $$ K_N\supset Q^{\alpha -\varepsilon_1}\leqno (5.10)$$
with probability greater than $1-2^{-n+1}$. Then, Theorem 3.9
shows that $$ K_N\supset F^{\alpha
-\varepsilon_1-\varepsilon_2}\cap\gamma C\leqno (5.11)$$ with
probability greater than $1-2^{-n+1}$, where $\varepsilon_2\leq 3
\log n/n$. Finally, Lemma 4.2 shows that there exists
$\varepsilon_3\leq 6/n$ such that $K_N$ satisfies $$ |\partial
(F^{\alpha +\varepsilon_3})\cap \gamma C\cap
K_N|\leq\tfrac{1}{2}|\partial (F^{\alpha +\varepsilon_3})\cap
\gamma C |\leqno (5.12)$$ with probability greater than
$1-10^{-2}$. Thus, there exists a 0/1  polytope $K_N$ which
satisfies all the above.

Next, apply Lemma 5.2 with $\beta =\alpha
-\varepsilon_1-\varepsilon_2$ and $\varepsilon
=\varepsilon_1+\varepsilon_2+\varepsilon_3$: If $A$ is a facet of
$K_N$ and $H_A$ is the corresponding halfspace (which has interior
disjoint from $K_N$), then $$ |\partial (F^{\alpha +\varepsilon_3
})\cap\gamma C\cap H_A|\leq \big (
3(\varepsilon_1+\varepsilon_2+\varepsilon_3) n\big )^{(n-1)/2}\,
|S^{n-1}|.\leqno (5.13)$$ It follows that
\begin{eqnarray*}
f_{n-1}(K_N)\, \big (3(\varepsilon_1+\varepsilon_2+\varepsilon_3)
n\big )^{(n-1)/2}\, |S^{n-1}| &\geq & \sum_A|\partial (F^{\alpha
+\varepsilon_3 })\cap\gamma C\cap H_A|\\ &=& \left |\big (\partial
(F^{\alpha +\varepsilon_3})\cap \gamma C\big )\setminus K_N\right
|\\ &\geq & \tfrac{1}{2}|\partial (F^{\alpha +\varepsilon_3})\cap
\gamma C |.
\end{eqnarray*}
Now apply Lemma 5.1 with $\beta =\alpha +\varepsilon_3$ to get $$
f_{n-1}(K_N)\, \big (3(\varepsilon_1+\varepsilon_2+\varepsilon_3)
n\big )^{(n-1)/2}\geq \tfrac{1}{2}\left (
(1-\gamma^2)\sqrt{2\alpha n}\right )^{n-1}.\leqno (5.14)$$ Since
$\alpha n=\log N$ and
$(\varepsilon_1+\varepsilon_2+\varepsilon_3)n\leq 12\log n$, this
shows that $$ f_{n-1}(K_N)\geq\left ( \frac{c_1(\gamma )\log
N}{\log n}\right )^{n/2}\leqno (5.15)$$ and the proof is complete.
$\hfill\Box $

\medskip

\noindent {\bf Proof of Theorem 1.1}: We can apply Theorem 1.2
with $N\geq \exp \big (c_1 n/\log n\big )$ where $c_1>0$ is an
absolute constant. This shows that there exists a 0/1 polytope $P$
in ${\mathbb R}^n$ with $$ f_{n-1}(P)\geq\left (
\frac{c_2n}{\log^2 n}\right )^{n/2},\leqno (5.16)$$ which is our
lower bound for $g(n)$. $\hfill\Box $

\bigskip

\bigskip

\noindent \textsc{D.\ Gatzouras}: Agricultural University of
Athens, Mathematics, Iera Odos 75, 118 55 Athens, Greece.
\textit{E-mail:} {\tt gatzoura@aua.gr}

\medskip

\noindent \textsc{A.\ Giannopoulos}: Department of Mathematics,
University of Crete, Heraklion 714 09, Crete, Greece.
\textit{E-mail:} {\tt giannop@fourier.math.uoc.gr}

\medskip

\noindent \textsc{N.\ Markoulakis}: Department of Mathematics,
University of Crete, Heraklion 714 09, Crete, Greece.
\textit{E-mail:} {\tt math2002@math.uoc.gr}

\end{document}